\documentclass[12pt]{amsart}

\newtheorem{theorem}{Theorem}
\newtheorem{lemma}{Lemma}
\newtheorem{proposition}{Proposition}

\newcommand{\ltwo}{L^2({\mathbb R})}
\newcommand{\bbz}{\mathbb{Z}}
\newcommand{\bbrn}{\mathbb{R}^n}
\newcommand{\Cal}{\mathcal}
\newcommand{\Bb}{\mathbb}
\newcommand{\comm}{\mathcal{C}_{\psi}(D,T)}

\begin{document}

% MAKE THE TITLE
\title{Applications of the Wavelet Multiplicity Function}
\author{Eric Weber}
\address{Department of Mathematics, University of Colorado,
Boulder, CO, 80309-0395} 
\email{webere@euclid.colorado.edu}
\thanks{This research was supported in part by a National Science Foundation grant, DMS9801658}
\begin{abstract}
This paper examines the wavelet multiplicity function.  An explicit formula for the multiplicity function is derived.  An application to operator interpolation is then presented.  We conclude with several remarks regarding the wavelet connectivity problem.
\end{abstract}
\maketitle

% BEGIN THE MAIN BODY OF THE PAPER

\section{Introduction} \label{S:Intro}
In this paper we are considering wavelets in the classic sense: a function $\psi \in \ltwo$ is a wavelet if the following is an orthonormal basis for $\ltwo$:
\[ 
\{2^{\frac{j}{2}} \psi (2^{j}x + n) : j,n \in \Bb{Z} \}.
\]
We define the operators $D$, $T$ on $\ltwo$ as: $Df(x) = \sqrt{2} f(2x)$ and $Tf(x) = f(x + 1)$.  Hence, a wavelet is a complete wandering vector for the unitary system $\{ D^n T^l : n,l \in \Bb{Z} \}$, see \cite{DL}.

Every wavelet can be associated with a \emph{Generalized Multiresolution Analysis}, or GMRA.  Indeed, define the subspaces $V_j = \{D^n T^l \psi : n < j,\ l \in \bbz\}$, then it is routine to verify that these subspaces satisfy the following four conditions:
\begin{enumerate}
\item $V_j \subset V_{j+1}$,
\item $DV_j = V_{j+1}$,
\item $\cap_{j \in \bbz} V_j = \{0\}$ and $\cup_{j \in \bbz} V_j$ has dense span in $\ltwo$,
\item $V_0$ is invariant under $T$. \label{I:4}
\end{enumerate}
In this paper, $V_0$ will often be referred to as the ``core space'' of $\psi$, as in \cite{BMM}.

Item~\ref{I:4} is of particular interest, since it yields a representation of the integers on $V_0$ via translations.  This representation is called the \emph{core representation}.  By Stone's theorem, there is a projection valued measure associated to the representation.  Then, that projection valued measure generates a probability measure $\mu$ and a \emph{multiplicity function} on the interval $[-\pi, \pi]$.  The probability measure and multiplicity function determine a unitary representation uniquely up to unitary equivalence.

Given two wavelets $\psi$ and $\eta$, each has an associated core space, which will be denote by $V_0^{\psi}$ and $V_0^{\eta}$, respectively.  Then, two wavelets $\psi$ and $\eta$ are said to be \emph{core equivalent} if their core representions are unitarily equivalent.  In \cite{BMM}, it is shown that the probability measure $\mu$ for any wavelet is absolutely continuous with respect to Lebesgue measure, which we shall denote by $\lambda$, hence the multiplicity function completely determines the representation.  Thus, two wavelets are core equivalent if and only if their associated multiplicity functions are the same.

For convenience of notation, if $f$ is measurable, for every $x \in \Bb{R}$, define a sequence on $\Bb{Z}$, denoted by $\vec{f}(x)$, by $\vec{f}(x)[k] = f(x + 2 k \pi)$.  If $f \in \ltwo$, then for almost every $x$, $\vec{f}(x) \in l^2(\Bb{Z})$.
  
\section{An Explicit Formula} \label{S:M=D}

In this section we will derive an explicit formula for the wavelet multiplicity function in terms of the wavelet itself.  We begin by decomposing the core representation into cyclic subrepresentations.  Cyclic representations have multiplicity functions that take on only the values $0$ and $1$, and the multiplicity function of a direct sum of representations is the sum of the individual multiplicity functions.  For $j > 0$, let $\psi_j = D^{-j}\psi$.  Let $g_1 = \psi_1$, and let $Y_1$ be the cyclic subspace generated by $g_1$.  Let $g_2 = P_{Y_1}^{\perp}\psi_2$, and let $Y_2$ be the cyclic subspace generated by $g_2$.  Recursively define $g_j$ to be the projection of $\psi_j$ onto the perpendicular complement of $\oplus_{n < j} Y_n$, and $Y_j$ to be the cyclic subspace generated by $g_j$.  Be definition, each $Y_j$ determines a cyclic subrepresentation, with cyclic vector $g_j$.

\begin{proposition}
\begin{enumerate}
\item $V_0 = \oplus_{j > 0} Y_j$,
\item $m = \sum_{j = 1}^{\infty} m_j$, where $m_j$ is the multiplicity function of the cyclic represention on $Y_j$.
\end{enumerate}
\end{proposition}
\begin{proof}
By definition, the $Y_j$'s are orthogonal and are subspaces of $V_0$.  Since the translations of the $\psi_j$'s spans $V_0$, it suffices to show that they are contained in this direct sum.  But note that $\psi_1$ is in $Y_1$, and then $\psi_2$ can be written as $g_2 + f_2$, where $f_2 \in Y_1$ since $g_2$ is obtained by a projection.  By the recursive definition of the $g_j$'s, we get that $\psi_j$ is in the direct sum, and item 1 is established.

Item 2 follows from the general fact that the muliplicity function for a representation is the sum of the multiplicity functions for orthogonal subrepresentations.
\end{proof}

Since $g_j$ is a cyclic vector, it generates a positive definite function, $p_j(l) = \langle T^l g_j, g_j \rangle$.  By Bochner's theorem, there exists a measure $\mu_j$ whose Fourier-Stieltjes transform is $p_j$.  Since $\mu$ is absolutely continuous with respect to Lebesgue measure, the measure $\mu_j$ is also.  Let $h_j$ be the Radon-Nikodym derivative of $\mu_j$ with respect to Lebesgue measure.  Since the subrepresentation on $Y_j$ is cyclic, $m_j$ takes on only the values 0 and 1; in fact, $m_j = \chi_{supp(h_j)}$.  Furthermore, since $m = \sum_{j = 1}^{\infty} m_j$ we have then that $m = \sum_{j=1}^{\infty} \chi_{supp(h_j)}$.

\begin{proposition}
Let $h_j$ be as above.  Then: 
\begin{equation}
h_j(\xi) = 2 \pi \|\Vec{\widehat{g_j}}(\xi) \|^2 \label{E:h_j}
\end{equation}
\end{proposition}
\begin{proof}
We have:
\begin{align}
\int_{0}^{2 \pi} e^{-in\xi} h_j(\xi) d\lambda
	&= \int_{0}^{2 \pi} e^{-in\xi} d\mu_j \notag \\
	&= \widehat{\mu_j}(-n) \notag \\
	&= p_j(-n) \notag \\
	&= \langle T^{-n} g_j, g_j \rangle \notag \\
	&= 2 \pi \int_{\Bb{R}} e^{-in\xi} \widehat{g}_j(\xi) \overline{\widehat{g}_j(\xi)} d\lambda \notag \\
	&= \int_{0}^{2 \pi} e^{-in\xi} 2 \pi \|\Vec{\widehat{g_j}}(\xi) \|^2 d\lambda \notag
\end{align}
\end{proof}

By definition, 
\begin{equation}
g_j(x) = \psi_j(x) - w_j(x) \label{E:gf}
\end{equation}
where $w_j$ is the unique element in $\oplus_{k < j} Y_k$ such that $g_j \perp \oplus_{k < j}Y_k$.  Additionally, since $w_j$ can be expressed in terms of the translates of the $g_k$'s, by taking the Fourier Transform of both sides of \ref{E:gf}, we get:
\[ \widehat{g_j}(\xi) = \hat{\psi}_j(\xi) - \sum_{k < j} \eta_{j,k}(\xi) \widehat{g_k}(\xi) \]
where the $\eta_{j,k}$ are $2 \pi$-periodic measurable functions.

\begin{proposition}
If $\eta_{j,k}$ is as above, then
\[ \eta_{j,k}(\xi) = \frac{\langle \Vec{\widehat{\psi_j}}(\xi), \Vec{\widehat{g_k}}(\xi) \rangle}{\| \Vec{\widehat{g_k}}(\xi) \|^2} \]
where this is interpreted to be 0 when the denominator is 0.
\end{proposition}
\begin{proof}
First notice that since $\eta_{j,k}\widehat{g_k} \in \ltwo$, $\eta_{j,k} (\xi)\| \Vec{\widehat{g_k}}(\xi) \|^2 \in L^1([0, 2 \pi])$.  Furthermore, $\eta_{j,k} \widehat{g_k} \in \widehat{Y_k}$, indeed, $\eta_{j,k} \widehat{g_k}$ is the function such that $\widehat{\psi_j} - \eta_{j,k} \widehat{g_k} \in \widehat{Y_k}^{\perp}$.  Hence, $\langle \eta_{j,k} \widehat{g_k}, e^{-in\cdot}\widehat{g_k} \rangle = \langle \widehat{\psi_j}, e^{-in\cdot}\widehat{g_k} \rangle$.  Hence,
\begin{align}
\int_{0}^{2 \pi} \eta_{j,k} (\xi) \| \Vec{\widehat{g_k}}(\xi) \|^2 e^{-in\xi} d\lambda
	&= \int_{\Bb{R}} \eta_{j,k} (\xi) \widehat{g_k}(\xi) \overline{\widehat{g_k}(\xi)} e^{-in\xi} d\lambda \notag \\
	&= \langle \eta_{j,k} \widehat{g_k}, e^{-in\cdot}\widehat{g_k} \rangle \notag \\
	&= \langle \widehat{\psi_j}, e^{-in\cdot} \widehat{g_k} \rangle \notag \\
	&= \int_{\Bb{R}} \widehat{\psi_j}(\xi) \overline{\widehat{g_k}(\xi)} e^{-in\xi} d\lambda \notag \\
	&= \int_{0}^{2 \pi} \langle \Vec{\widehat{\psi_j}}(\xi), \Vec{\widehat{g_k}}(\xi) \rangle e^{-in\xi} d\lambda \notag
\end{align}
\end{proof}

This gives us that
\begin{equation}
\widehat{g_j}(\xi) = \widehat{\psi_j}(\xi) - \sum_{k < j} \frac{ \langle \vec{\widehat{\psi_j}}(\xi), \vec{\widehat{g_k}}(\xi) \rangle}{\| \vec{\widehat{g_k}}(\xi)\|^2} \widehat{g_k} (\xi). \label{E:g_j}
\end{equation}
Since the inner product in equation \ref{E:g_j} is invariant under $2 \pi$ translations, we have:
\begin{equation}
\vec{\widehat{g_j}}(\xi) = \Vec{\widehat{\psi_j}}(\xi) - \sum_{k < j} \langle \Vec{\widehat{\psi_j}}(\xi), \Vec{u}_k (\xi) \rangle \vec{u}_k(\xi) \label{E:gbar}
\end{equation}
where
\[ \Vec{u_k} (\xi) = \frac{\vec{\widehat{g_k}}(\xi)}{\|\vec{\widehat{g_k}}(\xi)\|_2} \]
if the norm is non-zero.

As we have mentioned above, the multiplicity function is the sum of the multiplicity functions for each cyclic subspace $Y_j$, each of which is the characteristic function of the support of $h_j$.  Hence, by equation \ref{E:h_j}, the multiplicity function is precisely the number of non-zero sequences $\vec{\widehat{g_j}}(\xi)$.

Now, let us examine more closely equation~\ref{E:gbar}.  Note that $g_1 = \psi_1$, so $\vec{\widehat{g_1}}(\xi) = \vec{\widehat{\psi_1}}(\xi)$ for almost all $\xi$, and $u_1$ is the normalization of that vector.  Furthermore, $\vec{\widehat{g_2}}(\xi) = \Vec{\widehat{\psi_2}}(\xi) - \langle \Vec{\widehat{\psi_2}}(\xi), \Vec{u_1} (\xi) \rangle \vec{u_1}(\xi)$ which is the Graham-Schmidt orthogonalization of $\Vec{\widehat{\psi_1}}(\xi)$ and $\Vec{\widehat{\psi_2}}(\xi)$, with $\Vec{u}_1(\xi)$ and $\Vec{u}_2(\xi)$ being normalized.  By the recursive definition of the $g_j$'s, equation~\ref{E:gbar} is actually the Graham-Schmidt orthogonalization of the $\Vec{\widehat{\psi_j}}(\xi)$'s.  Hence, $h_j(\xi) = 0$ if and only if $\Vec{\widehat{\psi_j}}(\xi)$ is in the linear span of the previous $\Vec{\widehat{\psi_n}}(\xi)$'s.  Therefore, $m(\xi)$ is the number of linearly independent vectors in the collection $\{\vec{\widehat{\psi_j}}(\xi)\}$, i.e. the dimension of the subspace spanned by those vectors.

Let $\psi$ be a wavelet on $\ltwo$, and define $D_{\psi}: [-\pi, \pi] \to \Bb{R}$ as follows: 
\[
D_{\psi}(\xi) = \sum_{j=1}^{\infty} \sum_{k \in \bbz} |\hat{\psi}(2^j (\xi + 2 \pi k))|^2.
\]
This function is called the dimension function, as in \cite{HW}.

\begin{theorem} \label{T:M=D}
Let $\psi$ be a wavelet, and let $m : [-\pi, \pi] \rightarrow \Bb{Z^+}$ be its associated multiplicity function.  Then
\[
m(\xi) = D_{\psi}(\xi).
\]
\end{theorem}

\begin{proof}
It is shown in \cite{HW} that the dimension function is integer valued, and in fact is the dimension of a subspace of $l^2(\bbz)$.  Let $\vec{\widehat{\Psi_j}}(\xi)$ be a sequence on $\Bb{Z}$ given by $\vec{\widehat{\Psi_j}}(\xi)[k] = \widehat{\psi}(2^j (\xi + 2 \pi k))$.  It is also shown that $D_{\psi}(\xi)$ is the dimension of the subspace of $l^2(\bbz)$ spanned by $\{\vec{\widehat{\Psi_j}}(\xi) : j > 0 \}$.

We have that both the dimension function and the multiplicity function describe the dimension of some subspace of $l^2(\Bb{Z})$.  The spanning vectors are different; however, the subspaces are the same.  Indeed, $\vec{\widehat{\Psi_j}}(\xi)[k] = \widehat{\psi}(2^j (\xi + 2 \pi k))$ where as $\Vec{\widehat{\psi_j}}(\xi)[k] = \widehat{\psi_j}(\xi + 2 \pi k) = 2^{\frac{j}2} \hat{\psi}(2^j (\xi + 2 \pi k))$.  Hence, $\vec{\widehat{\Psi_j}}(\xi) = 2^{\frac{j}2} \vec{\hat{\psi_j}}(\xi)$, and so the spaces spanned by them are the same.
\end{proof}

\section{Operator Interpolation of Wavelets} \label{S:OIW}

In Dai and Larson, the notion of a local commutant is introduced.  The local commutant for the unitary system $\Cal{U}$ is define as:
\[
\Cal{C}_{\psi}(\Cal{U}) = \{V \in \Cal{B}(\ltwo): VU(\psi) - UV(\psi)=0 \ \forall \ U \in \Cal{U} \}
\]
We shall denote by $\comm$ the local commutant for the unitary system $\{ D^n T^l : n,l \in \Bb{Z} \}$.  It is shown in~\cite{DL} that $\eta$ is a wavelet if and only if there exists a unitary operator $U \in \comm$ such that $U(\psi) = \eta$.

\begin{lemma} \label{L:UV}
Let $\psi$ be a wavelet, let $U \in \comm$ be a unitary operator, let $\eta = U(\psi)$, and let $V$ be another operator such that $VU$ is also in $\comm$.  Then, $V$ is in $\Cal{C}_{\eta}(D,T)$.
\end{lemma}
\begin{proof}
We need to show that $VD^nT^l(\eta) = D^nT^lV(\eta)$ for all $n, l \in \bbz$.  We have:
\begin{align}
VD^nT^l(\eta) = VD^nT^lU(\psi)& = VUD^nT^l(\psi) \notag \\
	&= D^nT^lVU(\psi) = D^nT^lV(\eta). \notag
\end{align}
\end{proof}

\begin{theorem} \label{T:Ueq}
Let $U$ be a unitary operator such that $U^n \in \comm$ for all $n \in \bbz$.  Then, we have a sequence of wavelets, $U^n(\psi)$, and they are all core equivalent.
\end{theorem}

\begin{proof}
Each wavelet can be associated to a GMRA; let $\psi^{(i)}$ denote the wavelet $U^i(\psi)$, and let $V_0^{i}$ denote the $V_0$ core space for the wavelet $\psi^{(i)}$.  We need to construct an intertwining operator between $V_0$ and $V_0^{(i)}$.

Let $Y$ denote the closed subspace spanned by $\cup_{i \in \bbz} (V_0^{i})^\perp$.

\begin{lemma} \label{L:UT}
If $x \in Y$, then $U^i T^n(x) = T^n U^i(x)$.
\end{lemma}

\textsc{Proof of Lemma.}
It suffices to establish the lemma for a generating vector of $Y$.  Let $x \in \cup_{i \in \bbz} (V_0^{i})^\perp$, then $x \in (V_0^{i_0})^\perp$ for some $i_0$.  Let $\eta$ denote the wavelet $\psi^{(i_0)}$. We have: 
\[
x = \sum_{j=0}^{\infty} \sum_{k \in \bbz} \langle x, D^j T^k \eta \rangle D^j T^k \eta.
\]
By Lemma~\ref{L:UV}, $U^i \in \Cal{C}_\eta(D,T)$.  Hence:
\begin{align}
U^i T^n (x) &= U^i T^n \sum_{j=0}^{\infty} \sum_{k \in \bbz} \langle x, D^j T^k \eta \rangle D^j T^k \eta \notag \\
	&= U^i \sum_{j=0}^{\infty} \sum_{k \in \bbz} \langle x, D^j T^k \eta \rangle D^j T^{2^j n + k} \eta \notag \\
	&= \sum_{j=0}^{\infty} \sum_{k \in \bbz} \langle x, D^j T^k \eta \rangle D^j T^{2^j n + k} U^i \eta \notag \\
	&= T^n \sum_{j=0}^{\infty} \sum_{k \in \bbz} \langle x, D^j T^k \eta \rangle D^j T^k U^i \eta \notag \\
	&= T^n U^i \sum_{j=0}^{\infty} \sum_{k \in \bbz} \langle x, D^j T^k \eta \rangle D^j T^k  \eta \notag \\
	&= T^n U^i (x). \notag
\end{align}
as required.

Since $(V_0^{(i)})^{\perp} \subset Y$ for all $i$, we have $Y^{\perp} \subset V_0^{(i)}$.  Write $V_0 = Y^{\perp} \oplus \tilde{Y}$ and $V_0^{(i)} = Y^{\perp} \oplus \tilde{Y}^{(i)}$.  Define $S: V_0 \to V_0^{(i)}$ by $S = P_{Y^{\perp}} + U^i P_{\tilde{Y}}$.

A routine computation shows that $U(V_0) = V_0^{(1)}$.  Since $U \in \Cal{C}_{\psi^{(i)}} (D,T)$ for all $i$ the subspace $Y$ is invariant under $U$.  Indeed, $U(Y) = Y$.  As such, $Y^{\perp}$ is also invariant under $U$.  Likewise, both $Y$ and $Y^{\perp}$ are invariant under $T$.

Note that $U$ maps $V_0$ unitarily onto $V_0^{(i)}$, and leaves $Y^{\perp}$ fixed, hence $U$ maps $\tilde{Y}$ unitarily onto $\tilde{Y}^{(i)}$.  This proves that $S$ is a unitary operator.  Furthermore, as noted above, all of these subspaces in question are invariant under translation, hence by lemma~\ref{L:UT} shows that $S$ commutes with translations.  Therefore, $S$ is the required intertwining operator.
\end{proof}

Dai and Larson consider the case when a collection of wavelets $\{\psi_i: i \in I \}$ has the property that each unitary corresponding with each pair of wavelets has all of its powers in the local commutant of the first wavelet.  Additionally, if all wavelets gotten by the powers of these unitaries are back in this collection of wavelets, this collection is called an \emph{interpolation family of wavelets}.  Applying theorem~\ref{T:Ueq} to this collection of wavelets, we get the following theorem.

\begin{theorem} \label{T:Feq}
If $\Cal{F} = \{\psi_i : i \in I\}$ is an interpolation family of wavelets, then they are all core equivalent.
\end{theorem}

Additionally, Dai and Larson present \emph{operator interpolation} which can be used to generate new wavelets.  The process is as follows: consider a unitary operator in $\comm$ such that all of its powers are back in $\comm$.  Additionally, suppose that $U^k = I$ for some $k$.  Then the collection $\{U^j(\psi): j < k\}$ is an interpolation family of wavelets.  Finally, assume that $U$ normalizes $\{D,T\}'$ in the sense that $U^* \{D,T\}' U = \{D, T\}'$.  Then this interpolation family is said to \emph{admit operator interpolation}.  Then, if $A_j \in \{D,T\}'$, we can consider the operator 
\[ V = \sum_{j = 1}^k A_j U^j. \]  This operator is back in $\comm$, and if it is unitary, then $V\psi$ is again a wavelet.  This wavelet is called an \emph{operator interpolated wavelet}.  It is shown in~\cite{DL}, chapter 5, that $V$ has all of its powers in the local commutant.  Again applying theorem~\ref{T:Ueq}, we have the following theorem.

\begin{theorem} \label{T:OIW}
If $\psi$ is an operator interpolated wavelet, interpolated from the interpolation family $\Cal{F}$, then $\psi$ is core equivalent to all of the wavelets in $\Cal{F}$.
\end{theorem}

The following theorem is an easy consequence of theorem~\ref{T:Ueq}, and will be discussed further in section~\ref{S:WCP}.

\begin{theorem} \label{T:Aeq}
If $\mathcal{A} \subset \comm$ is a von Neumann algebra, then all wavelets that are parametrized by unitary operators in $\mathcal{A}$ are core equivalent.
\end{theorem}

\textsc{Remark 1.}  By theorem~\ref{T:M=D}, theorems~\ref{T:Ueq}, \ref{T:Feq}, \ref{T:OIW}, and \ref{T:Aeq} can be restated in terms of the dimension function.

\textsc{Remark 2.} Nowhere in this section have we explicitly or implicitly used the fact that these are functions on the real line.  Indeed, the analysis extends trivially to $\bbrn$.  Hence, theorems~\ref{T:Ueq}, \ref{T:Feq}, \ref{T:OIW}, and \ref{T:Aeq} are valid if the wavelets are taken to be in $L^2(\bbrn)$.  The proofs are exactly the same.

\section{The Wavelet Connectivity Problem} \label{S:WCP}

The wavelet connectivity problem is the question of whether the collection of all wavelets forms a path connected subset of the unit sphere of $\ltwo$ with respect to the norm topology.  Using the unitary operators in a von Neumann algebra is useful for solving this problem since the unitary group of a von Neumann algebra is path connected in the operator topology.  Hence, of particular interest is when the local commutant contains a von Neumann algebra.  In leiu of theorem~ \ref{T:Aeq}, a necessary condition for a von Neumann algebra to be in the local commutant is that all of the associated wavelets are core equivalent.  Hence, as far as the local commutant parametrizing connected subsets of wavelets, it can only parametrize those that are core equivalent.

A sufficient condition for a von Neumann algebra to be in the local commutant would be of great importance, as it may yield new path connected subsets of the collection of wavelets.  However, the ``converse'' of theorem \ref{T:Ueq} is false, i.e. if $\psi$ and $\eta$ are core equivalent, then the unitary in $\comm$ that maps $\psi$ to $\eta$ does not necessarily have all of its powers back in $\comm$.  We shall now provide an example.

If $W_1$ and $W_2$ are wavelet sets, then there exists a measurable bijection $\sigma: W_1 \to W_2$ that is effected by $2 \pi$ translations.  This $\sigma$ can be extended to a measurable bijection of the entire real line by defining $\sigma(x)$ to be $2^{-n}\sigma(2^n x)$, where $2^n x \in W_1$.  Define $U$ a unitary operator on $\ltwo$ by $Uf = f \circ \sigma^{-1}$.  It can be shown that $U \in \Cal{C}_{\hat{\psi}_{W_1}}(\hat{D}, \hat{T})$.  Additionally, $U^2 \in \Cal{C}_{\hat{\psi}_{W_1}}(\hat{D}, \hat{T})$ if and only if $\sigma^2|_{W_1}$ is again effected by $2 \pi$ translations.  For a complete discussion, see~\cite{DL}, chapters 4 and 5.

Let $W_1 = [-\frac{\pi}{4}, -\frac{\pi}{8}) \cup [\frac{15\pi}{8}, \frac{15\pi}{4})$, let $W_2 = -W_1$, and let $\psi_{W_1}$ and $\psi_{W_2}$ be their associated wavelets, respectively.  It is easy to verify that both $\psi_{W_1}$ and $\psi_{W_2}$ are associated to MRA's, hence they are core equivalent.

A routine calculation shows that:
\begin{equation}
\sigma(\xi) = \notag
  \begin{cases}
    \xi - 2 \pi, &\text{$\xi \in [-\frac{\pi}{4}, -\frac{\pi}{8}) \cup 
				[\frac{17\pi}{8}, \frac{18\pi}{8})$} \\
    \xi - 4 \pi, &\text{$\xi \in [\frac{15\pi}{8}, \frac{17\pi}{8})$} \\
    \xi - 6 \pi, &\text{$\xi \in [\frac{18\pi}{8}, \frac{30\pi}{8})$} \\
  \end{cases}
\end{equation}
Consider $A = [\frac{17\pi}{8}, \frac{273\pi}{128}) \subset W_1$.  Then, $A - 2 \pi = [\frac{\pi}{8}, \frac{17 \pi}{128}) \subset [\frac{\pi}{8}, \frac{\pi}{4}) \subset W_2$, hence, $\sigma(\xi)|_A = \xi - 2 \pi$.  Now let's calculate $\sigma^2(\xi)|_A = \sigma(\xi)|_{A - 2\pi}$.  Note that $16(A - 2\pi) = [2 \pi, \frac{17\pi}{8})$, so then if $\xi \in A$, $\sigma^2(\xi) = \sigma(\xi - 2\pi) = \frac{1}{16}\sigma(16(\xi -2 \pi)) = \frac{1}{16}[(16(\xi - 2\pi) - 4\pi] = \xi - 2\pi - \frac{\pi}{4}$.  Hence, $U^2$ given by $\sigma^2$ cannot be in $\Cal{C}_{\hat{\psi}_{W_1}}(\hat{D}, \hat{T})$.

\bigskip

ACKNOWLEDGEMENTS:  First, the author wishes to thank the organizers of this Special Session, Dr.\ Larry Baggett and Dr.\ David Larson, for the invitation to participate.  Secondly, it was brought to the author's attention during the Special Session that theorem~\ref{T:M=D} was proven independently by Zioma Rzeszotnik and Marcin Bownik, based on work by De Boor, DeVore, and Ron on shift invariant subspaces of $\ltwo$.  Finally, the author wishes to thank Qing Gu for the counterexample presented in section~\ref{S:WCP}.

% Here is the Bibliography

\end{document}